\documentclass[12pt,twoside]{amsart}
\usepackage{amsmath}
\usepackage{amssymb}
\usepackage{amscd}
\usepackage{xypic}
\newmuskip\pFqmuskip

\usepackage{amsthm}
\usepackage{graphicx}
\usepackage{enumerate}
\usepackage{graphicx}
\usepackage{tikz-cd}
\usepackage[new]{old-arrows}
\usepackage[cmtip,all]{xy}

\usepackage{amsmath, amsthm, amssymb}

\xyoption{all}
\title[Tropical Normal Functions - Higher Abel-Jacobi Invariants of Tropical cycles]{Tropical Normal Functions - Higher Abel-Jacobi Invariants of Tropical cycles}

\author{Mohammad Reza Rahmati}
\thanks{}
\address{ Universidad De La Salle Bajío, Campestre - Le\'on, Guanajuato, Mexico.
\hfill\break 
\hfill\break \\
\hfill\break }
\email{mrahmati@cimat.mx}

\newcommand{\comments}[1]{}

\keywords{Family of tropical manifolds, Tropical Hodge structure, Tropical normal function, Tropical intermediate Jacobian, Tropical Abel-Jacobi Invariants, Tropical Hodge bundle}

\subjclass{05E14, 05E40, 05E45}
\begin{document}

\begin{abstract}
We consider the variation of tropical Hodge structure (TVHS) associated to families of tropical varieties. The family of the tropical intermediate Jacobians of the associated tropical Hodge structure defines a bundle of tropical Jacobians, whose sections we call the tropical normal functions. We define formal sequential derivatives of these functions on the base with respect to the natural Gauss-Manin connection as the Hodge theoretic invariants detecting tropical cycles in the fibers. The associated invariants which are defined inductively are the higher Abel-Jacobi invariants in the tropical category. They naturally identify the tropical Bloch-Beilinson filtration on the tropical Chow group. We examine this construction on the moduli of tropical curves with marked points, in order to study the tropical tautological classes in the tautological ring of $\mathcal{M}_{g,n}^{\text{trop}}$. The expectation is the nontriviality of these cycles could be examined with less complexity in the tropical category. The construction is compatible with the tropicalization functor on the category of schemes, and the aforementioned procedure will also provide an alternative way to examine the relations in the tautological ring of $\mathcal{M}_{g,n}$ in the schemes category. 
\end{abstract}

\maketitle

\section*{Summary}

The study of the tautological classes in the Chow ring of the moduli space of curves of genus $g$ with $n$-marked points has been under investigation in the last decays in different fields of researches in mathematics and also in Physical sciences. The combinatorial framework of these objects makes them interesting for the experts in different areas. On the other hand the difficulties of doing specific calculations with the algebraic and homological cycles in the Chow ring of $\mathcal{M}_{g,n}$ have forced the experts to search for new techniques that enable them, how much more characterize the classes. This article is a text toward this purpose to provide new methodology to study the generators and relations in the tautological ring of $\mathcal{M}_{g,n}$. 

The motivation we follow toward the above problem is through the tropical category. In fact the tropical manifolds and their tropical cycles can be constructed through the tropicalization functor applied to algebraic varieties. The idea is, the non-triviality of algebraic cycles in the Chow group of $\mathcal{M}_{g,n}$ can be checked out via this functor. In other words if the cycle has non-trivial invariants after the tropicalization procedure then the original algebraic cycle must be non-trivial. A natural expectation is tropicalizing the objects and the invariants may make the calculations simpler.

Historically Normal functions were introduced by H. Poincare. Higher cycle classes and Abel-Jacobi invariants were used by Green and Griffiths to show the non-triviality of a special zero cycle presented by Faber-Pandharipande. The cycle was originally from the pull back of a tautological cycle over the self product of the universal curve $C_g^2$. They calculate the infinitesimal invariants of the normal function associated to a thickening of the cycle in a generic family. The resulting cycle is homologicaly trivial and its Abel-Jacobi image is zero. But the second infinitesimal invariant of the thickening cycle is non-zero, which shows it is rationally non-trivial. 

The purpose of this text is to introduce tropical normal functions as sections of intermediate Jacobians of tropical variation of Hodge structure, TVHS. The work of this paper maybe explained in two ways. The first, one of the current research programs in tropical geometry involves possible ways to generalize geometric concepts to tropical category. In this way our definition of tropical normal function introduces a new concept in tropical geometry, as an extension the concept of normal functions in the category of complex VHS on quasi projective schemes over $\mathbb{C}$ into the tropical category. We also define higher normal functions and higher Abel-Jacobi maps in the tropical sense. The second, is that we use this concept as a tool to study algebraic cycles both in the tropical schemes in its own and also in the original category of schemes over $\mathbb{C}$. This not only gives a method to investigate intersection theory and Chow groups of tropical varieties, but also a way to check out relations and non-triviality of the original cycles in the scheme category. By the way using the properties of the tropicalization functor; the relations in the category of schemes imply same relations between the corresponding tropical cycles. We apply this simple fact to study tropical tautological classes in the tropical Chow group of the moduli of tropical curves with marked points. 

Let $\mathbb{T}=\mathbb{R} \cup -\infty$ be the tropical semifield. Tropical varieties are finite dimensional polyhedral complexes with specific affine structure. For example $\mathbb{T}\mathbb{P}^n$ is a smooth projective tropical variety homeomorphic to $n$-simplex. The restriction of the tropical structure to the interior of a $k$-dimensional face of $\mathbb{T}\mathbb{P}^n$ makes it into $\mathbb{R}^k$. Projective tropical varieties are polyhedral complexes in $\mathbb{T}\mathbb{P}^n$. Examples of this happens as tropical limit of an algebraic family $X_t \subset \mathbb{C}\mathbb{P}^n, \ t \in \mathbb{C}$ when $|t| \to \infty$. The fundamental theorem of tropical geometry states that $\log_{|t|}(X_t) \subset \mathbb{T}\mathbb{P}^n$ converges to an $n$-dimensional weighted polyhedral complex. A smooth projective tropical variety is a closed subset of $\mathbb{T}\mathbb{P}^n$ which on standard affine charts restricts to a smooth polyhedral complex in $\mathbb{T}^n$. [\cite{MS, GS, MZ, FS, K, RGST}].

We shall work with two dual lattices $M$ and $N=\text{Hom}(M,\mathbb{Z})$. Let $T=\text{Spec}(\mathbb{C}[M])$, be the algebraic torus. Let $Y \subset T$ be a closed subset. The tropicalization of $Y$ is defined via the map

\begin{equation}
    \begin{aligned}
    &\text{trop}:T^{\text{an}} \longrightarrow N_{\mathbb{R}}\\
    &\langle \text{trop}(x), \alpha \rangle =\log |\varpi^{\alpha}|_x
    \end{aligned}
\end{equation}

\noindent 
where $\varpi$ is the character in $K[M]$ corresponding to $m \in M$ and $|\ .\ |_x:K[M] \longrightarrow \mathbb{R}$ is the semi-norm extending $|.|$ on $K=\overline{\mathbb{C}(t)}$. The tropicalization is a non-archimedean moment map. [\cite{UL1}, \cite{UL2}]. 

The singular (co)homology of tropical varieties is defined in a natural way as for topological spaces. There is a natural generalization (definition) of Dolbeault cohomology to the tropical category (see Section 2 below). 
If $\mathfrak{X} \subset \mathbb{C}\mathbb{P}^n \times S^*$ is a holomorphic 1-parameter family of schemes over a disc $S^*$, with a tropical limit $X_{\infty} \subset \mathbb{T}\mathbb{P}^n$, then one can define 

\begin{equation}
H_{\text{trop}}^{p,q}= Gr_{2p}^WH^{p+q}(X_{\infty}, \mathbb{Q})    
\end{equation}

\noindent 
where $H^k(X_{\infty}, \mathbb{Q})$ is the limit mixed Hodge structure. We have the Steenbrink-Illusie spectral sequence with first cohomology 

\begin{equation}
E_1^{l,k-l}=\bigoplus_{s \geq \text{max}\{0,l\}}H^{k+l-2s}(\mathfrak{X}^{(2s-l)}, \mathbb{Q})[l-s])    
\end{equation}

\noindent 
that converges to $H_{trop}^{p,q}$, where $X^{(k)}$ is the $k$-skeleton of $X$, [\cite{IKGM}]. A similar formula holds over the complex Hodge structure. This shows that the Hodge filtration is well behaved under tropicalization functor.
 
Applying the tropicalization functor to a family $\mathfrak{X} \to S$ yields a family of tropical varieties defined over a base tropical variety $S$. The (tropical) cohomologies of the fiber varieties $X_s$ vary continuously with $s$, and form a vector bundle with fibers

\begin{equation}
    H_{\text{trop}}^k(X_s, \mathbb{R}):=\bigoplus_{p+q=k}H_{\text{trop}}^{p,q}(X_s, \mathbb{R})
\end{equation}

\noindent 
where $H_{\text{trop}}^{p,q}(X_s, \mathbb{R})$ are tropical Dolbeault cohomologies calculated via tropical Dolbeault complex, [\cite{IKGM}]. A Hodge type filtration can be defined as

\begin{equation}
    F_pH_{\text{trop}}^k(X_s, \mathbb{R})= \bigoplus_{l \geq p}H_{\text{trop}}^{l,k-l}(X_s,\mathbb{R})
\end{equation}

\noindent 
The tropical intermediate Jacobians of the fiber varieties is defined by

\begin{equation}
  J_{\text{trop}}^k(X_s)=\frac{H_{\text{trop}}^{2k-1} }{F^k+H_{\text{trop}}^{2k+1}(X_s, \mathbb{Z})} 
\end{equation}

\noindent 
They constitute a fiber bundle over $S$ with torus fibers. We call the sections of this bundle, the tropical normal functions. We purpose to associate homological invariants $\Psi_i^{\text{trop}}(Z_{\text{trop}})$ to tropical cycles which define a filtration $L^{\bullet}$ on the tropical Chow groups $CH_{\text{trop}}^k$ of length not bigger than $k+1$. Specifically

\begin{equation}
    Z_{\text{trop}} \in L^i CH_{\text{trop}}^k(\mathfrak{X}) \ \Leftrightarrow \ \Psi_0^{\text{trop}}(Z_{\text{trop}})=\Psi_1^{\text{trop}}(Z_{\text{trop}})=...=\Psi_{i-1}^{\text{trop}}(Z_{\text{trop}})=0
\end{equation}

\noindent 
The invariant $\Psi_i^{\text{trop}}(Z_{\text{trop}})$ is defined when all the former ones are zero. These invariants are a sum of higher cycles with higher Abel-Jacobi classes. One investigate if these invariants are non-trivial for tropical cycles, starting from some index on. 

Let $\Gamma=(V,E)$ be a weighted graph, with weight function $w$ (see Sections 4). A marked tropical curve is a continuous map $f:(\Gamma, x_1,...,x_n) \to M_{\mathbb{R}}$ satisfying certain balancing conditions. The moduli of tropical curves of genus $g$ with $n$ marked points is connected to $\mathcal{M}_{g,n}$ via a tropicalization map 

\begin{equation}
trop:\mathcal{M}_{g,n} \longrightarrow \mathcal{M}_{g,n}^{trop}
\end{equation}

\noindent 
The map is basically studied in [\cite{UL1, UL2}, see also the refs. there]. We have the universal tropical curve $C_g^{\text{trop}} \to \mathcal{M}_g^{\text{trop}}$ defined similarly. By the Leray-Deligne spectral sequence  

\begin{equation}
H^k(C_{g,\text{trop}}^n,\mathbb{Q})=\bigoplus_{p+q=k}H^q(\mathcal{M}_g^{\text{trop}},R^pf_*\mathbb{Q})    
\end{equation}

\noindent
Using a Schur-Weyl duality argument the decomposition tells that the cohomologies $H^k(C_{g,\text{trop}}^n,\mathbb{Q})$ and the collections $H^q(\mathcal{M}_g^{\text{trop}}, V_{\lambda}^{\text{trop}})$ where $V_{\lambda}$ are local systems on $\mathcal{M}_g$ corresponding to highest weight representations of $Sp_{2g}\mathbb{C}$, have the same information. 
As a consequence the set of invariants gets decomposed via the duality

\begin{equation} 
\Psi_i^{\text{trop}}(Z_{\text{trop}})= \bigoplus_{\lambda} \Psi_i^{\text{trop}}(Z_{\text{trop}})^{\lambda}
\end{equation} 

\noindent 
The way that the invariants $\Psi_i^{\text{trop}}(Z_{\text{trop}})$ were defined tells that for a nontrivial cycle $Z^{\text{trop}}$, in the above decomposition starting from some point the invariants show to be nontrivial, which depends to representation theoretic properties of $Sp_{2g}\mathbb{C}$. That is, in this case we have two sorts of filtrations on the tropical Chow groups, one is the tropical Bloch-Beilinson filtration characterized by the Hodge theoretic invariants and the other from the Schur-Weyl duality. The interaction between these two filtrations is the motivation of this text.

\section{Overview of Tropical Varieties}
\label{sec1}

The materials in this section are well known and can be found in various texts [\cite{GS, MS, MR1, SP, EH, GK, MI2, K, RGST, AR, CA1, CA2, EKL, GA, GU1, GG, KS, PA1, PA2}]. We work over the tropical ring $\mathbb{R}^{\text{trop}}[x_1,...,x_n]$ where $\mathbb{R}^{\text{trop}}=\mathbb{R}$ with tropical operations. A tropical polynomial is of the form 

\begin{equation}
    f(x_1,...,x_n)=\sum_{i_1,...,i_n \in S}a_{i_1,...,i_n}x_1^{i_1}...x_{n}^{i_n}=\text{min} \{a_{i_1,...,i_n}+ \sum_l i_l x_l | (i_1,...,i_n \in S\}
\end{equation}

\noindent 
where $S$ is a finite set of exponents. The zero set $V(f) \subset \mathbb{R}^n$ is a tropical hypersurface. We employ the two latices 
\begin{equation} 
M=\mathbb{Z}^n, \qquad N =Hom(M,\mathbb{Z})
\end{equation} 
Denote $M_{\mathbb{R}}=M \otimes \mathbb{R}$ and $N_{\mathbb{R}}=N \otimes \mathbb{R}$ and by $\langle . , . \rangle $ the natural pairing between them. A tropical function is a map 
\begin{equation} 
f:M_{\mathbb{R}} \to \mathbb{R}, \qquad  \text{given by (11) where} \ S \subset N
\end{equation} 
The Newton polytope of $f$ is the convex hull 
\begin{equation} 
\Delta_S=\text{Conv}(S) \subset N_{\mathbb{R}}
\end{equation} 
Recall, an affine linear automorphism of $M_{\mathbb{R}}$ is given as $m \longmapsto Am+b , \ A \in GL_n(\mathbb{R}),\ b \in M_{\mathbb{R}}$. A tropical manifold is a real topological manifold $X$ with an atlas of coordinate charts $\phi_i:U_i \to M_{\mathbb{R}}$ such that the transitions $\phi_i \circ \phi_j^{-1}$ are affine linear automorphisms. We allow the manifolds to have boundaries and singularities. 

A polyhedron $\sigma \in \mathbb{R}^n$ is an intersection of finite closed half spaces. We denote its boundary with $\partial \sigma$ and its interior by $\text{int}(\sigma)$. We call a compact polyhedron a polytope. A polyhedral decomposition of a polyhedron $\Delta \subset N_{\mathbb{R}}$ is a set $P$ of polyhedra in $N_{\mathbb{R}}$ namely cells such that 

\begin{itemize}
\item $\Delta=\bigcup_{\sigma} \sigma$.
\item If $\sigma \in P$ and $\tau \in \sigma$ then $\tau \in P$.
\item $\sigma_1, \sigma_2 \in P$ then $\sigma \cap \sigma_2 \in P$. 
\end{itemize}

\noindent 
If the half spaces are defined over $\mathbb{Q}$ we call the polyhedral decomposition of lattice type. We denote by $P^{(i)}$ the collection of $i$-dimensional cells in $P$. If $\tau \in P$ define the open star of $\tau$ by 

\begin{equation}
    \text{star}(\tau)=\bigcup_{\sigma \supset \tau} \text{Int}(\sigma)
\end{equation}

\noindent 
A fan structure along $\tau$ is a continuous map 

\begin{equation}
St_{\tau}:\text{star}(\tau) \to \mathbb{R}^k, \qquad k=\dim X -\dim \tau 
\end{equation}

\noindent 
such that 

\begin{itemize}
    \item $St_{\tau}^{-1}(0)=\text{Int}(\tau)$.
    \item If $\tau \subset \sigma$ then $St_{\tau}|_{\text{Int}}(\sigma)$ is an affine submerssion onto its image defined over $\mathbb{Q}$, where $\dim(St_{\tau}(\sigma))=\dim(\sigma)-\dim(\tau)$.
    \item For $\tau \subset \sigma$ define $C_{\tau, \sigma}$ to be the cone over $St_{\tau}(\sigma \cap star(\tau))$. Then 
    
    \begin{equation}
        \Sigma_{\tau}=\{C_{\tau, \sigma}|\tau \subset \sigma \in P\} 
    \end{equation}
    
    \noindent 
    is a fan with $|\Sigma_{\tau}|$ convex. 
\end{itemize}

\noindent 
Two fan structures are equivalent if one could be transformed to the other by an automorphism in $GL_k(\mathbb{Z})$. If $\tau \subset \sigma$ then a fan structure along $\tau$ induces a fan structure along $\sigma$ via the composition 

\begin{equation}
\text{star}(\sigma) \hookrightarrow \text{star}(\tau) \to \mathbb{R}^k \to \mathbb{R}^k/L_{\sigma} \cong \mathbb{R}^l 
\end{equation}

\noindent 
where $L_{\sigma} \subset \mathbb{R}^k$ is the linear span of $C_{\tau,\sigma}$. We may apply this definition along vertices of the polyhedral complex $X$ with a plyhedral decomposition $\Delta$. Let $v \in P$ be a vertex of $X$, and $U_v \subset \text{star}(v)$ be an open neighborhood such that for any codimension 1 cone $\rho$ the intersection $U_v \cap \rho$ is connected, then 

\begin{equation}
\{ \text{Int}(\sigma)|\sigma \in P^{\text{max}}\} \  \bigcup \ \ \{U_v|v \in  P^{(0)}\} \end{equation}

\noindent 
is an open cover of $X \setminus \Delta$. We define an affine structure on $X \setminus \Delta$ via

\begin{equation}
\begin{aligned}
&\psi_{\sigma}:int(\sigma) \to M_{\mathbb{R}}, \qquad \sigma \in P^{\text{max}}\\
&\psi_v: U_v \hookrightarrow \text{star}(v) \to \mathbb{R}^{\dim X}
\end{aligned}
\end{equation}

\noindent 
Say a collection of fan structures $\{St_v|v \in P^{[0]}\}$ are compatible if the induced fan structures on $\tau$ from all $ v \in \tau$ are equivalent. A tropical manifold is a pair $(X, P)$ where $X$ is a tropical affine manifold with singularities obtained from polyhedral decomposition $P$ of $X$ and compatible collection $\{St_v| v \in P^{(0)}\}$ of fan structures.  

We can also work out tropical varieties via their coordinate functions. The tropical projective space is defined as 

\begin{equation}
\begin{aligned}
   & \mathbb{T}\mathbb{P}^n = \{(x_0,...,x_n)| x_i \in \mathbb{T}, \text{not all} \ x_j =\infty \}/ \cong \\
   & (x_0,...,x_n) \cong (y_0,...,y_n) \ \ \Leftrightarrow \ \ x_j=\lambda  +y_j, \ \text{some} \ \lambda \in \mathbb{T}^*=\mathbb{R}
\end{aligned}    
\end{equation}

\noindent 
Let $K=\overline{\mathbb{C}(t)}$ the algebraic closure of rational functions in the indeterminate $t$ (field of Puiseux series in $t$). Assume 
\begin{equation} 
\text{val}:(K \setminus 0)^n \longrightarrow \mathbb{Q}^n
\end{equation} 

\noindent 
be a valuation map, where we may assume is the order map at $0$, i.e it assigns the least power of $t$ that appear in a Puiseux series. In the category of varieties over $K$, to any ideal $I$ of $K[x_1^{\pm 1},...,x_n^{\pm 1}]$, we can associate the zero set $V(I) \subset (K \setminus 0)^n$. The image of $V(I)$ under $\text{val}$ is a subset of $\mathbb{Q}^n$. By taking its tropical closure the resulting subset of $\mathbb{R}^n$ is a tropical variety $\mathbb{T}(I)$. A tropical algebraic variety is any subset of $\mathbb{R}^n$ of the form 

\begin{equation}
    \mathbb{T}(I)=\overline{\text{val}\left (V(I) \right )}, \qquad I \leq K[x_1^{\pm 1},...,x_n^{\pm 1}]]
\end{equation}

\noindent 
An ideal $I$ is homogeneous if all the monomials appearing in a generator of $I$ have the same total degree. Then $I$ defines a variety in $\mathbb{P}_K^{n-1} \setminus \{ x_j=0 , \forall j\}$. Its image under the map $\text{val}$ is a subset of the tropical projective space $\mathbb{T}\mathbb{P}^{n-1}$. A tropical projective variety is a subset of $\mathbb{T}\mathbb{P}^{n-1}$ of the form 

\begin{equation}
\mathbb{T}(I)=\overline{\text{val} \left (V(I) \right )}   
\end{equation}

\noindent 
where $I$ is a homogeneous ideal of $K[x_1^{\pm 1},...,x_n^{\pm 1}]$. A tropical curve is an embedded graph in $\mathbb{T}\mathbb{P}^2$ which is dual to the regular subdivision $\Delta$ of support of a tropical polynomial $f$. 

We can construct tropical varieties via the functor called tropicalization applied to the schemes in algebraic geometry. The functor is an analogue of the moment map in toric geometry. Let $X=X(\Delta)$ to be a toric variety, with $\Delta$ is a rational polyhedral fan in $N_{\mathbb{R}}$. The tropicalization on $\Delta$ is given by

\begin{equation}
\text{trop}_{\Delta}:U_{\sigma}^{\text{an}} \longrightarrow (N_{\sigma})_{\mathbb{R}}   \end{equation}

\noindent 
where $U_{\sigma}=\text{Spec} K[S_{\sigma}], \ S_{\sigma}=\sigma^{\vee} \cap M$ and $N_{\sigma}=Hom(S_{\sigma},\mathbb{R})$. The compactification of $N_{\sigma}$ is given by $Hom(S_{\sigma},\overline{\mathbb{R}})$. The image of tropicalization map is a partial compactification $N_{\mathbb{R}}(\Delta)$ of $N_{\mathbb{R}}$ determined by $\Delta$

\begin{equation}
\begin{aligned}
    &\text{trop}_{\Delta}:U_{\sigma}^{\text{an}} \longrightarrow N_{\mathbb{R}}(\sigma)\\
    &\text{trop}_{\Delta}(x)(s)=-\log|\varpi^s|_x
    \end{aligned}
\end{equation}

\noindent
(where $\varpi, |.|$ was introduced in the introduction) is a non-archemidean analytic moment map.The suffix 'an' is the analytification, [\cite{UL1, UL2}].

\section{Tropical Homology and Algebraic Cycles}

The singular homology of a tropical variety $X$, is defined naturally as usual varieties [\cite{IKGM}]. The tropical Dolbeault cohomology groups are defined as the cohomology of Dolbeault complex $(A^{\bullet, \bullet}, d', d'')$ of smooth real differential forms on the analytic space $X^{an}$, where

\begin{equation}
A_{X^{\text{an}}}^{p,q}(U)=C^{\infty}(U) \otimes \bigwedge^pT{X^{\text{an}}} \otimes \bigwedge^q T^*{X^{\text{an}}} 
\end{equation}

\noindent 
The operators $d', d''$ are formal diffentiations on the corresponding coordinate variables on $\bigwedge^pT{X^{\text{an}}}$ and $\bigwedge^q T^*{X^{\text{an}}}$, respectively. Specifically 

\begin{equation}
    H^{p,q}(X^{an})=\frac{\ker \left ( d'': A^{p,q} \to A^{p,q+1} \right )}{d''(A^{p,q-1})}
\end{equation}

\noindent 
The cohomology of the tropical variety is defined by 

\begin{equation} 
H^k(X^{\text{an}}):=\bigoplus_{p+q=k}H^{p,q}, \qquad H^{q,p}=(H^{p,q})^{\vee}
\end{equation} 

\noindent 
There is a well defined map namely tropical monodromy 

\begin{equation}
    N_{\text{trop}}:H_{\text{trop}}^{p,q}(X)_{\mathbb{R}} \to H_{\text{trop}}^{p-1,q+1}(X)_{\mathbb{R}}
\end{equation}

\noindent 
which is induced from a map at the level sheaves 

\begin{equation}
\begin{aligned}
&A^{p,q}(X^{\text{an}}) \to  A^{p-1,q+1}(X^{\text{an}}) \\
\omega = \sum_{|I|=p,|J|=q}&a_{I,J}d'x_{i_1} \wedge ... \wedge d'x_{i_p} \wedge  d''x_{j_1} \wedge ... \wedge d''x_{j_q} \longmapsto \\
& \sum_{l=1}^p \sum_{I,J}(-1)^{p-l}a_{I,J}d'x_{i_1} \wedge .. \wedge \widehat{d'x_{i_l}} \wedge . . \wedge d'x_{i_p} \wedge d''x_{i_l} \wedge d''x_{j_1} \wedge ... \wedge d''x_{j_q}
\end{aligned}
\end{equation}

In the rest of the section we briefly introduce the notion of tropical algebraic cycles and tropical Chow group, [\cite{MR1, SP, SH, FR, AC, AM1, BA1, BN, MI2, MZ, FS, K, AR, GG}]. We first introduce the concept of the weight function on polyhedrans. Assume $X$ is a tropical variety. Let $\Delta^{(k)}$ be the set of all cones of codimension $k$ on $X$. If $\sigma \in \Delta^{(k)}, \ \tau \in \Delta^{(k+1)}$, Let $N_{\sigma}$ be the lattice span of $\sigma$, and let $n_{\sigma,\tau} \in \sigma$ be the primitive integer vector whose image generates the 1-dimensional lattice $N_{\sigma}/N_{\tau}$. A weight function in codimension $k$ is a map 

\begin{equation}
    w:\Delta^{(k)} \longrightarrow \mathbb{Z}
\end{equation}

\noindent 
such that for any $\tau \in \Delta^{(k+1)}$ and any vector in $u \in \tau^{\perp} \cap \mathbb{Z}^n$, 

\begin{equation}
\sum_{\text{codim}_{\sigma}(\tau)=1} w(\sigma) \langle u , n_{\sigma, \tau} \rangle =0  
\end{equation}

\noindent 
A polyhedral complex $\sigma \in \Delta^{(k)}$ in a tropical variety $X$ is called an algebraic cycle if the following balancing condition are satisfied

\begin{equation}
\sum_{\text{codim}_{\sigma}(\tau)=1}w(\sigma)\ n_{\sigma,\tau}=0    
\end{equation}

\noindent 
We denote the free abelian group on all such $\sigma \in \Delta^{(k)}$ by $Z_k^{\text{trop}}(X)$. A tropical regular function on a tropical cycle $Z \in Z_k^{\text{trop}}(X)$ is a piece-wise affine function with integral slopes which is the restriction of a concave function on $\mathbb{R}^n$, i.e of the form 

\begin{equation}
f=\text{min}_{j \in J}\ \{\langle \alpha_j,. \rangle +\beta_j\}, \qquad \alpha_j \in \mathbb{Z}^n, \ \beta_j \in \mathbb{Q}    
\end{equation}

\noindent 
Then the divisor of $f$ is 

\begin{equation}
\text{div}(f)=\left [ \left ( \bigcup_{i=1}^{k-1} Z^{(i)}, w_f \right ) \right ] \in Z_{k-1}^{\text{trop}}(X)   
\end{equation}

\noindent 
where 

\begin{equation}
\begin{aligned}
    &w_f: X^{(k-1)} \longrightarrow \mathbb{Z}\\
    &\tau \longmapsto \sum_{\sigma \supset \tau}f_{\sigma} \big (w(\sigma)n_{\sigma,\tau} \big )-f_{\tau}\left (\sum_{\sigma \supset \tau}w(\sigma) n_{\sigma,\tau} \right )
\end{aligned}    
\end{equation}

\noindent
where $f_{\sigma}:N_{\sigma} \to \mathbb{R}$ denotes the linear part of the affine function $f|_{\sigma}$ and $\text{codim}_{\sigma}\tau =1$. We define $R_k$ the subgroup of divisors $(f)$ where $f$ is a function on $Z \in Z_{k}$. Finally the tropical Chow group of degree $k$ is defined as $CH_k=Z_k/R_k$. The total tropical Chow ring is 

\begin{equation} 
CH_{\text{trop}}^{\bullet} (X)=\bigoplus_k \ CH_k(X)
\end{equation} 

\noindent 
We shall consider the convention $CH_k=CH^{n-k}$. We have 

\begin{itemize}
    \item Tropical cycle class map 
    
    \begin{equation}
    \begin{aligned}
       & CH_{\text{trop}}^p({X}) \stackrel{\cong}{\longrightarrow}H_{\text{trop}}^{p,p}({X})\\
       & CH_{\text{trop}}({X}) \cong \bigoplus_p H_{\text{trop}}^{p,p}({X})
    \end{aligned}    
    \end{equation}
    
    \item Tropical intersection pairing 
    \begin{equation}
    \begin{aligned}
     &CH_{\text{trop}}^k ({X}) \times  CH_{\text{trop}}^{n-k}({X})  \longrightarrow CH_{\text{trop}}^n({X}) \stackrel{\deg}{\longrightarrow} \mathbb{Z} \\
     &(a,b) \longmapsto a.b \longmapsto \deg(a.b)
     \end{aligned}
    \end{equation}
    \item Tropical Hard Lefchetz
    
    \begin{equation}
    \begin{aligned}
        &CH_{\text{trop}}^k({X}) \stackrel{\cong}{\longrightarrow} CH_{\text{trop}}^{n-k}({X})\\
        &a \longmapsto l^{n-2k}.a
    \end{aligned}    
    \end{equation}
    
    \noindent 
    and the pairing 
    
    \begin{equation}
    \langle .,l^{n-k} \rangle : P_k^{\text{trop}}(X) \times P_k^{\text{trop}}(X) \longrightarrow \mathbb{Z}    
    \end{equation}
    
    \noindent 
    is nondegenerate, where $l$ is the class of a hyperplane, and 
    
    \begin{equation}
        P_k^{\text{trop}}: =\ker \left ( l^{n-k+1}:CH_{\text{trop}}^k(X) \longrightarrow CH_{\text{trop}}^{n-k}(X) \right )
    \end{equation}
\end{itemize}

\noindent
One may formulate the above properties in terms of the Hodge classes $H_{\text{trop}}^{p,p}$ instead. 

We explain the case of curves as the simplest case [\cite{GG, GA, EKL, AR, RGST, VI, MZ, MI2, HM, CA1, CA2, CA3, AM1, AC}]. Let $C=(V,E)$ be a graph. If $v \in V$ is a vertex its valance is the number of edges connected to $v$. A metric graph is a graph with length function

\begin{equation}
    l:E \longrightarrow \mathbb{R}_{>0}
\end{equation}

\noindent 
The first betti number of $C$ is called its genus, which is

\begin{equation}
g(C)=|E|-|V|+|\pi_0(C)|    
\end{equation}

\noindent 
A tropical curve is a connected metric graph $C$ where $val(v) \geq 2$ for all $v \in V$. A divisor on $C$ is an element of the free abelian group generated by the points of $|C|$, i.e $D=\sum a_iP_i$. A rational function on $|C|$ is a continuous piece-wise linear function $f:|C| \to \mathbb{R}$ with integer slopes. One can adjust this by specifying slopes on the graph edges. The orver at $p$ of $f$ is defined as 

\begin{equation}
    \text{ord}_p(f)=\sum_i\frac{\partial f}{\partial x_i}(p)
\end{equation}

\noindent 
The divisor of $f$ is given by $(f)=\sum_P\text{ord}_p(f)P$. The degree of the divisor $(f)$ is $0$. The tangent space $T_pC$ at $p$ is the set of derivations $\partial/\partial x_i$ of $C$ at $p$ corresponding to the tangent directions at $p$. A holomorphic differential form on $C$ is a collection 

\begin{equation}
w_p:T_pC \longrightarrow \mathbb{R} \qquad \text{such that} \ \sum_iw_p(\frac{\partial}{\partial x_i})=0 , \ \ (p \in C)   
\end{equation}

\noindent 
By choosing $g$ points $P_1,...,P_g$ of valance 2 in $C$ together with a primitive tangent vector $\partial/\partial x_i$ at each of these points we can define an assignment 

\begin{equation}
    w \longmapsto \left (w_{P_1}(\frac{\partial}{\partial x_1}), ..., w_{P_1}(\frac{\partial}{\partial x_g}) \right )
\end{equation}

\noindent 
For any path $\gamma \in C$ and $w \in \Omega^1(C)$ we can define the integral $\int_{\gamma} w$  by pulling back the 1-form on an interval. we have the embedding 

\begin{equation}
\begin{aligned}
    &H_1^{\text{trop}}(C, \mathbb{Z}) \hookrightarrow \Omega_{\text{trop}}^1(C)\\ 
    &\gamma \longmapsto \left ( \int_{\gamma}: (w \longmapsto \int_{\gamma}w) \right )
    \end{aligned}
\end{equation}

\noindent 
We define the Jacobian of $C$ by 

\begin{equation} 
J^{\text{trop}}(C)=\frac{\Omega_{\text{trop}}^1(C)^*}{H_1^{\text{trop}}(C, \mathbb{Z})}=\frac{\mathbb{R}^g}{\mathbb{Z}^g}
\end{equation} 

\noindent 
via the above embedding.


\section{Tropical Normal Functions and (higher) Abel-Jacobi Invariants}

In this section we define tropical Hodge theoretic invariants detecting tropical cycles in the tropical Chow group of tropical varieties. The corresponding invariants in the category of quasiprojective schemes over $\mathbb{C}$ are already well known [\cite{KLM, K1, K2, K3, GG1, GG2, L, V1, VO1, VO2, WU, ZU1, ZU2, H, CGG}].
Assume $\mathfrak{X} \subset \mathbb{C}\mathbb{P}^n \times S^*$ is a holomorphic 1-parameter family of schemes over a disc $S^*$, with a tropical limit $X_{\infty} \subset \mathbb{T}\mathbb{P}^n$, then one can define 

\begin{equation}
H_{\text{trop}}^{p,q}= Gr_{2p}^WH^{p+q}(X_{\infty}, \mathbb{Q})    
\end{equation}

\noindent 
where $H^k(X_{\infty}, \mathbb{Q})$ is the limit mixed Hodge structure. We have the Steenbrink-Illusie spectral sequence with first cohomology 

\begin{equation}
E_1^{l,k-l}=\bigoplus_{s \geq \text{max}\{0,l\}}H^{k+l-2s}(\mathfrak{X}^{(2s-l)}, \mathbb{Q})[l-s])    
\end{equation}

\noindent 
that converges to $H_{trop}^{p,q}$, where $X_{\infty}^{(k)}$ is the $k$-skeleton of $X_{\infty}$, [\cite{IKGM}]. A similar formula holds over the complex Hodge structure. This shows that the Hodge filtration is well behaved under tropicalization.  

We consider a family of projective tropical varieties defined over a quasi-projective tropical variety $S$. Because the discussion is local such a family can be simply be a product $X \times S$. The family can be defined by an affine map with parameters varying in a tropical variety. The fiber varieties $X_s$ vary continuously with $s$, and their cohomologies form a vector bundle with fibers

\begin{equation}
    H_{\text{trop}}^k(X_s, \mathbb{R})=\bigoplus_{p+q=k}H_{\text{trop}}^{p,q}(X_s, \mathbb{R}), \qquad H_{\text{trop}}^{q,p}=\left ( H_{\text{trop}}^{p,q} \right )^{\vee}
\end{equation}

\noindent 
We call the decomposition the tropical Hodge structure (THS). Also analogous to the ordinary category of quasi-projective schemes over $\mathbb{C}$ we can define TMHS when similar decomposition is satisfied on each graded pieces of the weight filtration. Tropical mixed Hodge structures TMHS form an abelian category. This category has $Hom$ and the first $Ext^1$ functor. The tropical Hodge filtration is defined similar to the complex case

\begin{equation}
    F_pH_{\text{trop}}^k(X_s, \mathbb{R})= \bigoplus_{l \geq p}H_{\text{trop}}^{l,k-l}(X_s,\mathbb{R})
\end{equation}

\noindent 
The family of $\mathbb{R}$-vector spaces $H_{\text{trop}}^k(X_s, \mathbb{R})$ as $s$ varies constitute a local system of $\mathbb{R}$-vector spaces $\mathcal{H}_{\text{trop}}$ over the topological space $S$, which is defined over $\mathbb{Q}$. By the Riemann-Hilbert correspondence this local system defines a flat connection 

\begin{equation} 
\nabla :\mathcal{H}_{\text{trop}} \to \mathcal{H}_{\text{trop}} \otimes \Omega_S^{\text{trop}, 1}
\end{equation} 

\noindent 
whose kernel is $\mathcal{H}_{\text{trop}}$. The interaction between this connection and the intermediate jacobians of $\mathcal{H}_{\text{trop}}$ is what we purpose to discuss. 

The tropical intermediate Jacobians is defined by

\begin{equation}
  J_{\text{trop}}^k(X_s)=\frac{H_{\text{trop}}^{2k-1} }{F^k+H_{\text{trop}}^{2k+1}(X_s, \mathbb{Z})} 
\end{equation}

\noindent 
Similar formulas holds over the complex Hodge structure. In other words by what we said above
\begin{equation}
J_{\text{trop}}^k(X_s)=\text{trop} (J^k(X_s))
\end{equation}
The formula of J. Carlson on the relation between $J^k$ and extensions of Hodge structures in complex Hodge theory translate to  
\begin{equation} 
J_{\text{trop}}^k(X_s)=\text{Ext}_{\text{TMHS}}^1(\mathbb{Z}, H_{\text{trop}}^{2k-1})
\end{equation} 
The family of $J_{\text{trop}}^k(X_s)$ constitute a bundle of toruses over $S$ which we denote by $\mathcal{J}_{\text{trop}}^k$. We call the sections of this bundle the tropical normal functions. 

The tropical Abel-Jacobi maps are defined similar to the case of curves [ see Section 2 above], as 

\begin{equation}
\begin{aligned}
    &CH_{\text{trop}, hom}^k(X) \stackrel{AJ^{\text{trop}}}{\longrightarrow}J^k(X)=\frac{(F^k H_{\text{trop}}^{2k-1}(X, \mathbb{C})^{\vee}}{\int_{H_{2k-1}^{\text{trop}}(X, \mathbb{Z})}(.)}\\
    & Z=\partial \gamma \ \longmapsto \int_{\gamma}\  (.)
\end{aligned}    
\end{equation}

\noindent 
The suffix $hom$ means homologically trivial cycles. We can define a sequence of invariants for the tropical cycles which define a filtration $L^{\bullet}$ on the tropical Chow groups $CH_{\text{trop}}^k$ of length not bigger than $k+1$. Specifically

\begin{equation}
    Z_{\text{trop}} \in L^i CH_{\text{trop}}^k(\mathfrak{X}) \ \Leftrightarrow \ \Psi_0^{\text{trop}}(Z_{\text{trop}})=\Psi_1^{\text{trop}}(Z_{\text{trop}})=...=\Psi_{i-1}^{\text{trop}}(Z_{\text{trop}})=0
\end{equation}

\noindent 
In fact the invariant $\Psi_i^{\text{trop}}(Z_{\text{trop}})$ is defined when all the former ones are zero. We can consider A Leray filtration $L^i$ whose graded parts are as follows 

\begin{equation}
\begin{aligned}
Gr_L^iHg_{\text{trop}}^k(\mathfrak{X})=Hom_{\text{TMHS}}(\mathbb{Q}(-k),H^i(S,R^{2k-i}\pi_*\mathbb{Q}))\\
Gr_L^iJ_{\text{trop}}^k(\mathfrak{X})=\frac{Ext_{\text{TMHS}}(\mathbb{Q}(-k),H^{i-1}(S,R^{2k-i}\pi_*\mathbb{Q}))}{Hom_{\text{TMHS}}(\mathbb{Q}(-k),H_{\text{trop}}^{2k-1}(\mathfrak{X},\mathbb{Q}))}
\end{aligned}
\end{equation}

\noindent
The above two pieces can be collapsed to define a filtration on tropical Chow group of $\mathfrak{X}$, where its $i$-th graded piece defines the $i$-th cycle class map on the tropical cycles, which we denote by $cl_{\text{trop}}^i=[.]_i^{\text{trop}}$. In fact these invariants are tied with higher Abel-Jacobi invariants $[AJ(.)]_{i-1}$ with a shift according to the second grading. We record this issue in a split short exact sequence of Abelian groups which maps these cycles as 

\begin{equation}
   [AJ(.)]_{i-1}^{\text{trop}} \longmapsto \Psi_i(.)^{\text{trop}} \longmapsto [.]_i^{\text{trop}} 
\end{equation}

\noindent
We have the following 

\begin{itemize}
    \item $[.]_i^{\text{trop}}$ is the image of $[Z]^{\text{trop}}$ under the projection
    
    \begin{equation}
     Hg^p (H_{\text{trop}}^{2p} (\mathfrak{X})) \to Hg^p(Gr_L^i H_{\text{trop}}^{2p} (\mathfrak{X})) \to Gr_L^iHg_{\text{trop}}^p(\mathfrak{X})
    \end{equation}
    
    \item $[AJ(.)]_{i-1}^{\text{trop}}$ is the image of $AJ^{\text{trop}}(.)$ under  
    
    \begin{equation}
     J^p (H_{\text{trop}}^{2p-1} ) \to J^p(Gr_L^i H_{\text{trop}}^{2p-1} (\mathfrak{X})) \to Gr_L^{i-1}J_{\text{trop}}^p(\mathfrak{X})
    \end{equation}
\end{itemize}

\noindent 
$[Z]_0^{\text{trop}}$ identifies $[Z]^{\text{trop}}$. If $[Z]_0^{\text{trop}}=[Z]_1^{\text{trop}}=0$ then $[AJ(.)]_0^{\text{trop}}$ is $AJ(.)^{\text{trop}}$. 

The aforementioned construction is justified by the tropicalization functor. The reader may know that the construction is well established in the category of schemes, \cite{K}. If we apply the tropicalization functor to a family of projective varieties over a quasi-projective base $\mathfrak{X} \to S$ we obtain a family of tropical varieties $trop(\mathfrak{X}) \to trop(S)$ where all the construction applies to. In fact the definition shows that the invariants just defined are the tropicalization of the corresponding invariant in the category of projective schemes. An example we may identify the BB-filtration mentioned above for a product of tropical curves. Let $X=C_1^{\text{trop}} \times ... \times C_n^{\text{trop}}, \ \sigma \in S_n$ and let 

\begin{equation} 
\sigma:(C_1^{\text{trop}} \times C_1^{\text{trop}} ) \times ...(C_1^{\text{trop}} \times C_1^{\text{trop}} ) \to X \times X
\end{equation}

\noindent 
the obvious map. Let 

\begin{equation}
    \pi_{\sigma}=\sum_{\alpha_1+...+\alpha_n=n}\pi_{C_1,\alpha_1} \times ... \times \pi_{C_n,\alpha_n}, \qquad \alpha_i =0, 1, 2
\end{equation}

\noindent 
where 

\begin{equation}
\pi_{C,2}=C^{\text{trop}} \times 0 , \qquad \pi_{C,0}=0 \times C^{\text{trop}}, \qquad \pi_{C,1}=\Delta-\pi_{C,2}-\pi_{C,0}    
\end{equation}

\noindent 
Then we have 

\begin{equation}
L_iCH_{\text{trop}}^n(X)=\bigcap_{\sigma \in \Xi_{i-1}^n}\ker (\pi_{\sigma})_*    
\end{equation}

\noindent
where $\Xi_{i-1}^n$ is the set of increasing functions $\epsilon:\{1,...,i\} \to \{1,...,n\}$. The projections onto the Chow-Kunneth components can be written as 

\begin{equation}
pr_{\lambda}=\sum_{{\stackrel {j < \lambda}{\sigma \in S_j}}} (-1)^{\lambda -j}(i_{\sigma} \circ \pi_{\sigma})
\end{equation}

\noindent
where $\pi_{\sigma}:X \to X_{\sigma}=C_{\sigma(1)}^{\text{trop}} \times ... \times C_{\sigma(i)}^{\text{trop}}$ and $i_{\sigma}:X_{\sigma} \hookrightarrow X$ are the obvious maps.

Let $\mathfrak{X}^{\text{trop}} \longrightarrow S^{\text{trop}}$ be a family in tropical category with fibers $X_s^{\text{trop}}$, and $Z^{\text{trop}}$ a tropical cycle of codimension $k$ such that $Z_s^{\text{trop}}=Z^{\text{trop}}.X_s^{\text{trop}}$ are also of codimension $k$ that is homologous to $0$ in tropical homology, then we can define a tropical normal function associated to $Z^{\text{trop}}$ by

\begin{equation}
 \nu_{Z^{\text{trop}}}(s)=AJ_{X_s}^{\text{trop}}(Z_s^{\text{trop}})  \in J_{\text{trop}}^k(X_s^{\text{trop}})=F_s^k \setminus H_s^{k, {\text{trop}}}/ H_{\mathbb{Z},s}^{k, {\text{trop}}}  
\end{equation}

\noindent 
We can lift the value of $\nu_{Z^{\text{trop}}}$ to a $\widetilde{\nu}_{Z^\text{trop}}$ in $H_s^{k, {\text{trop}}}$ locally. Then we have 

\begin{equation} 
\nabla (\widetilde{\nu}_{Z^{\text{trop}}})(s) \in F_s^{k-1}
\end{equation}

\noindent 
Thus we may define an infinitesimal invariant for $\nu_{Z^{\text{trop}}}$ by

\begin{equation}
    \delta ( \nu_{Z^{\text{trop}}})=\nabla (  \widetilde{\nu}_{Z^{\text{trop}}})
\end{equation}

\noindent 
The invariant we are interested in is its class as 

\begin{equation}
    \delta (\nu_{Z^{\text{trop}}}) \in H_{\text{trop}}^{k,k-1}/H_{\text{trop}}^{k-1,k-2}
\end{equation}

\noindent 
We can proceed to define inductively normal type functions ${\nu}_{Z^{\text{trop}}}^{(i), } \in [J_{\text{trop}}^k(X_s^{\text{trop}})]_i$ and the invariants 

\begin{equation}
\delta^{(i)}(\nu_{Z^{\text{trop}}}) \in H^i(H_{\text{trop}}^{k-i,k+i-1} , \nabla)  
\end{equation}

\noindent 
and so on. This how the above sequence of invariants works out. If at some stage one of these invariants turns to be non-zero, it follows that the original cycle must be non trivial.

\section{Tautological Classes in Moduli of Tropical Curves with Marked Points}

We apply the technology of the previous section to tautological classes in the tropical Chow ring of $\mathcal{M}_{g,n}^{\text{trop}}$, the moduli of tropical curves of genus $g$ with $n$ marked points, [\cite{MI1, IKE, MR, GS, FR, CA2, CHMR, CGP, GKM, MI3, RA, UL1, UL2}]. The idea is a family of tropical curves over a base can also be studied over the moduli of tropical curves, as a universal curve. Lets recall the definition of a tropical curve with marked points. We always assume a graphs $\Gamma=(V,E)$, is weighted, i.e. come with a weight function 

\begin{equation} 
w:E \longrightarrow \mathbb{N}
\end{equation}

\noindent 
A marked parametrized curve is a continuous map $f:(\Gamma, x_1,...,x_n) \to M_{\mathbb{R}}$ such that

\begin{itemize}
    \item If $w(E)=0$, then $f|_E$ is constant. Otherwise $f|_E$ is a proper embedding of $E$ into a line of rational slope in $M_{\mathbb{R}}$.
    \item If $v \in V$ and $E_1,...,E_k$ are the edges adjacent to $v$ and $e_i \in M$ be a primitive tangent vector to $f|_{E_i}$ then
    
    \begin{equation}
     \sum_i\ w(E_i)\ e_i=0   
    \end{equation}
\end{itemize}

\noindent 
The moduli of tropical curves of genus $g$ with $n$ marked points is connected to the ordinary $\mathcal{M}_{g,n}$ via a tropicalization map 

\begin{equation}
trop_{g,n}:\mathcal{M}_{g,n} \to \mathcal{M}_{g,n}^{trop}
\end{equation}

\noindent 
The map is studied in [\cite{UL1, UL2}]. We have the universal tropical curve $C_g^{\text{trop}} \to \mathcal{M}_g^{\text{trop}}$ where is the tropicalization of $C_g$. Take the $n$-fold fibration

\begin{equation} 
C_{g,\text{trop}}^n=(C_g^{\text{trop}})^n= C_g^{\text{trop}} \times_{\mathcal{M}_g^{\text{trop}}} ... \times_{\mathcal{M}_g^{\text{trop}}} C_g^{\text{trop}} \longrightarrow \mathcal{M}_g^{\text{trop}}
\end{equation} 

\noindent 
That is the moduli of tropical curves of genus $g$ with $n$ ordered marked points. By the Leray-Deligne spectral sequence of the fibration $f:(C_g^{\text{trop}})^n \to \mathcal{M}_g^{\text{trop}}$, we obtain a decomposition 

\begin{equation}
H^k(C_{g,\text{trop}}^n,\mathbb{Q})=\bigoplus_{p+q=k}H^q(\mathcal{M}_g^{\text{trop}},R^pf_*\mathbb{Q})    
\end{equation}

\noindent 
The above decomposition is the base of some analysis of cycle classes in the Chow ring of $C_{g,\text{trop}}^n$. It is more convenient to consider the tautological subrings of tropical cycles or their corresponding Hodge classes generated by the diagonals $\Delta_{ij}^{\text{trop}}$ of two marking points and the tropical classes $\psi_j^{\text{trop}}$ that are the chern classes of $n$ cotangent bundles at the marked points, and the Morita-Mumford-Miller classes $\kappa_d$. According to the major result in [\cite{PTY}] a Schur-Weyl duality argument applies with an application of the theorem of Ancona, [\cite{A}] to obtain

\begin{equation}
    CH_{\text{trop}}(C_{g,\text{trop}}^n)= \bigoplus_{|\lambda| \leq n}V_{\lambda}^{\text{trop}} \otimes \mathbb{L}_{\text{trop}}^{n_i}
\end{equation}

\noindent 
which descends to cohomology and also to tautological subalgebras. Using the decomposition in (80) the higher Abel-Jacobi invariants for the tropical cycles in the Chow group of $C_{g,\text{trop}}^n$ splits as in the following decomposition

\begin{equation}
    \begin{aligned}
    &L^iCH_{\text{trop}}^n(C_{g,\text{trop}}^n) \stackrel{\oplus (\pi_{\lambda})_*}{\longrightarrow} \bigoplus_{\mid \lambda \mid \leq n} L^iCH_{\text{trop}}^{n-n_{\lambda}}(S, \textbf{V}_{\lambda})\\
    &Gr_L^iH_{\text{trop}}^{2n}(C_{g,\text{trop}}^n,\mathbb{C})  \stackrel{\cong}{\longrightarrow} \bigoplus_{\mid \lambda \mid \leq n}Gr_L^iH_{\text{trop}}^{2n-2n_{\lambda}}(S,\textbf{V}_{\lambda})
    \end{aligned} 
\end{equation}

\noindent
Moreover, one has 

\begin{equation}
L^iCH_{\text{trop}}^n(C_{g,\text{trop}}^n)=\ker\{\bigoplus_{\mid \lambda \mid <i}(\pi_{\lambda})_*\}
\end{equation}

\noindent
If $\mathfrak{Z} \in CH_{\text{trop}}^n(C_{g,\text{trop}}^n)$, then

\begin{itemize}
\item $\Psi_i(\mathfrak{Z})^{\text{trop}}= \bigoplus_{\lambda} \Psi_i(\mathfrak{Z})^{{\text{trop}},\lambda}$
\item $[\mathfrak{Z}]_i^{\text{trop}} =\bigoplus_{\lambda} [\mathfrak{Z}]_i^{\text{trop},\lambda}$
\item $[AJ(\mathfrak{Z})]_{i-1}^{\text{trop}} =\bigoplus_{\lambda} [AJ(\mathfrak{Z})]_{i-1}^{\text{trop},\lambda}$
\end{itemize} 

\noindent 
We have 

\begin{equation}
\bigoplus[AJ(\mathfrak{Z})]_{i-1}^{\text{trop},\lambda} \hookrightarrow \Psi_i^{\text{trop}}(\mathfrak{Z})= \bigoplus_{\lambda} \Psi_i(\mathfrak{Z})^{{\text{trop}},\lambda}\twoheadrightarrow \bigoplus_{\lambda} [\mathfrak{Z}]_i^{\text{trop},\lambda}
\end{equation}

\noindent 
We consider the class 

\begin{equation}
FP_1^{\text{trop}}:=\pi_1^{\times 2}(\Delta_{12} . \psi_1)_{\text{trop}} \in CH_{\text{trop}}^2(C_g^2)
\end{equation}

\noindent
called the tropical Faber-Pandharipande cycle. It is a codimension 2 cycle in $C_{g, {\text{trop}}}^2$. In the scheme category the restriction of the $FP_1$ is a special cycle worked by Green and Griffiths in [\cite{GG1}]. If $C^{\text{trop}}$ is a general curve of genus $g \geq 4$. Let $\imath^{\Delta}: C^{\text{trop}} \hookrightarrow C^{\text{trop}} \times C^{\text{trop}}$. The ordinary tropical Faber-Pandharipande cycle 

\begin{equation}
Z_K^{\text{trop}}:=K_C^{\text{trop}} \times K_C^{\text{trop}} -(2g-2)\imath_*^{\Delta}K_C^{\text{trop}} \in CH_{\text{trop}}^2(C^{\text{trop}} \times C^{\text{trop}})
\end{equation}

\noindent
is $\stackrel{\text{rat}}{\ne} 0$, where $K_C^{\text{trop}}$ is the tropical canonical class of $C^{\text{trop}}$. For $g=2,3$ it is known to be rationally equivalent to zero, (see the ref.). It is homologically and also AJ-equivalent to zero, that is $[Z_K^{\text{trop}}]=Alb^{\text{trop}}(Z_K^{\text{trop}})=0$ for all $g$. Our construction together with the work of Green-Griffiths gives the following, 

\begin{itemize}
\item $[Z_K]_0^{\text{trop}}=[Z_K^{\text{trop}}]=0$
\item $[Z_K]_1^{\text{trop}}=0$ 
\item $[AJ(Z_K^{\text{trop}})]_0^{\text{trop}}:=AJ_{C \times C}^{\text{trop}}(Z_K)=0$
\item $Z_K^{\text{trop}} \in  L^2CH_{\text{trop}}^2(C^{\text{trop} , \times 2})$
\item $[Z_K]_2^{\text{trop}} =0$
\item $[AJ(Z_K^{\text{trop}})]_1^{\text{trop}} \ne 0$
\end{itemize}

\noindent
We give some other examples obtained from the application of tropicalization functor to the corresponding cycles in the scheme category. Define the $n^{\text{th}}$ tropical Faber-Pandharipande cycle 

\begin{equation}
FP_n^{\text{trop}} =\pi_1^{\times 2} (\Delta_{12}^n. \psi_1) \in CH_{\text{trop}}^{n+1}(C_{g, {\text{trop}}}^2)
\end{equation}

\noindent
The cycle $FP_n^{\text{trop}}$ has codimension $n+1$ in $C_{g, {\text{trop}}}^2$. A related cycle is

\begin{equation} 
\mathfrak{GS}^{\text{trop}}=\pi_1^{\times 3} \Delta_{123} \in CH_{\text{trop}}^2(C_{g, {\text{trop}}}^3)
\end{equation} 

\noindent 
called the tropical Gross-Schoen cycle, and originally studied in [\cite{GSc}] in the scheme category. B. Gross and C. Schoen study the cycle 

\begin{equation}
Y=\Delta_{123}-\Delta_{12}-\Delta_{13}-\Delta_{23}+\Delta_1+\Delta_2+\Delta_3
\end{equation}

\noindent
where $\Delta_I=\{(x_1,...,x_n): x_i=x_j \ \text{if}\ i,j \in I \ \text{and} \ x_i=o \ \text{if} \ i \notin I\}$. The tropical diagonals are also defined analogously. Let $\mathfrak{Y}$ be the spread of $Y$ over $S$.  

\begin{equation}
(r_S^{[3]})^* \mathfrak{GS}^{\text{trop}}- \mathfrak{Y}^{\text{trop}} =\text{sum of $FP_1^{\text{trop}}$-cycles}
\end{equation} 

\noindent
where $(r_S^{[3]})$ is inclusion. It follows that

\begin{itemize}
\item $[\mathfrak{Y}]_0^{\text{trop}}=0$. 
\item $[AJ(\mathfrak{Y}^{\text{trop}})]_0^{\text{trop}} \ne 0$.
\item $[AJ(\mathfrak{Y}^{\text{trop}})]_1^{\text{trop}} \ne 0, \ \  (g \geq 3)$.
\end{itemize} 

\noindent
The application of the tropical functor with the result of [\cite{GSc}] predicts that when $C^{\text{trop}}$ is a general tropical curve of genus $g \geq 3$, then $\mathfrak{Y}^{\text{trop}} \stackrel{\text{rat}}{\ne} 0$. A more general example could be

\begin{equation}
\mathfrak{FP}(n,m)^{\text{trop}}:=\pi_1^{\times n} (\Delta_{12...n}.\psi_1^m)_{\text{trop}} \in CH_{\text{trop}}^{n+m-1}(C_{g, {\text{trop}}}^n)
\end{equation}

\noindent
The tropical cycle in 4.16 is the image of the cycle $\pi_1^{\times 3} (\Delta_{123})^{\times (n+2m-2)}$ under (tropicalized) correspondence map constructed in [\cite{PTY}, Lemma 12.4, page 43], via the projectors in the Brauer algebra such that

\begin{equation}
H_{\text{trop}}^{n+2m-2}(M_g, V^{\otimes 3(n+2m-2)}) \to H_{\text{trop}}^{n+2m-2}(M_g^{\text{trop}},V^{\otimes n})
\end{equation} 

\noindent
The compatibility of Abel-Jacobi invariants under pull back implies

\begin{equation}
[AJ(\mathfrak{Y}_S^{{\text{trop}}, \times (n+m-1)})]_{n+m-1}^{\text{trop}} \ne 0 \ \ \ \Rightarrow \ \ [AJ(\mathfrak{FP}(n,m)_S)]_{n+m-1}^{\text{trop}} \ne 0
\end{equation}
  
\noindent
which proves non triviality of $\mathfrak{FP}(n,m)_S^{\text{trop}}$ in appropriate generic case. 

One may consider universal normal functions on the moduli of tropical curves, [\cite{H}]. Let $\mathcal{A}_g^{\text{trop}}=\text{trop}(\mathcal{A}_g)$ be the moduli of  polarized tropical abelian varieties of dimension $g$. Let $f:\mathfrak{X} \to \mathcal{A}_g^{\text{trop}}$ be the universal tropical abelian variety. The local system $H=R^1f_*\mathbb{Z}$ is a variation of tropical Hodge structure of weight 1. We denote the pull back of $H$ under the period mapping $\mathcal{M}_{g,n}^{\text{trop}} \to \mathcal{A}_g^{\text{trop}}$ by the same symbol $H$. The family of intermediate Jacobians $J^{\text{trop}}(H)$ is isomorphic to $\mathcal{X}$. Another interesting TVHS is $\bigwedge^3 H$. There is an embedding of $H \hookrightarrow \bigwedge^3 H$ name,y $x \hookrightarrow x \wedge Q$ where $Q:H \otimes H \to \mathbb{Z}$ is the Poincare pairing. Denote $V=\bigwedge^3H/H$. 

We explain two sort of natural normal functions on $\mathcal{M}_{g,n}^{\text{trop}}$. The first is 
\begin{equation}
\begin{aligned}
    K_j&:\mathcal{M}_{g,n}^{\text{trop}} \longrightarrow J(H), \\
    K_j&([C^{\text{trop}},x_1,...,x_n])=(2g-2)x_j-K_C^{\text{trop}}, \ n \geq 1
\end{aligned}    
\end{equation}
The second is given by 
\begin{equation}
\begin{aligned}
    D_{j,k}&:\mathcal{M}_{g,n}^{\text{trop}} \longrightarrow J(H), \\
    D_{j,k}&([C^{\text{trop}},x_1,...,x_n])=x_j-x_k, \ n \geq 1
\end{aligned}    
\end{equation}
A conjectural question is if the group of normal functions of $J(H) \to \mathcal{M}_{g,n}^{\text{trop}}$ is freely generated by the above two series of normal functions. In fact this holds in the scheme category for the aforementioned moduli spaces. 

We can extend the aforementioned question for the TVHS of the tensor power $H^{\otimes n}$. From the decomposition 
\begin{equation}
    H_{\text{trop}}(C_{g,\text{trop}}^n)= \bigoplus_{|\lambda| \leq n}V_{\lambda}^{\text{trop}} \otimes \mathbb{L}_{\text{trop}}^{n_i}
\end{equation}
and additivity of $\text{Ext}$ it suffices we determine the set of normal functions for the THS $V_{\lambda}^{\text{trop}} \otimes \mathbb{L}_{\text{trop}}^{n_i}$. The conjecture is that this set is zero unless for the highest weights $\lambda=[1], [1,1,1]$ of $SP_g$-representations. 

We mention another remark following [\cite{PTY}] for the algebraic version, that a kind of formulas which descend from the scheme category to the tropical moduli is the formalism of FZ relations on the relations in the tautological ring of $\mathcal{M}_{g,n}^{\text{trop}}$, due to Janda-Pixton-Pandharipande-Zvonkine. That is the expression 
\begin{equation}
    \left[ \exp(-\{\log(A)\}_n\sum_{|P|=n}\prod_{S \in P}\{C_{|S|}\}_{\Delta_S} \right]_{z^r}
\end{equation}
vanishes in $CH_{\text{trop}}^r(C_{g, \text{trop}}^n)$. The expression reads as follows. $P$ runs over the partitions of $n$. The expression $C_n$ is defined inductively via the two power series 
\begin{equation}
A(z)=\sum_j \frac{6j)!}{(2j)!(3j)!}z^j, \qquad B(z)=\sum_j \frac{6j)!}{(2j)!(3j)!}\frac{6j+1}{6j-1}z^j   
\end{equation}
where $C_0=\log (A), \ C_1=B/A, \ C_{n+1}=(12z^2\frac{d}{dz}-4nz)C_n$. If $f(z)=\sum_ja_jz^j$ then $\{f\}_k=\sum_jk_ja_jz^j$ and $\{f\}_{\Delta S}=\sum_j(-1)^{|S|-1}\Delta_S\psi_S^{j-|S|+1}a_jz^j$, and $S \subset \{1,...,n\}$. Finding all the relations in the tautological ring $R_{\text{trop}}^r(C_{g, \text{trop}}^n)$ is a challenging research in this area.
\vspace{0.3cm}

\subsection*{Acknowledgments.} I learned about tropical Hodge bundle in the conference ICM 2018 Satellite: Tropical geometry and moduli spaces, on 05 Dec 2017 13 aug 2018 - 17 aug 2018 Cabo Frio, Rio de Janeiro and found about the above motivation in some conversations in the conference. I wanted to thank of the organizers of the conference. 

\subsection*{Contributions:} The normal functions associated to variation of tropical Hodge structure has been introduced together with the higher tropical cycle classes and Higher Abel-Jacobi Invariants of tropical cycles in the fibered variety. We have used these invariants to study generators and relations in the tropical tautological ring of the moduli or tropical curves of genus $g$ and $n$-marked points.

\end{document}